\documentclass[12pt]{amsart}

\usepackage{enumitem}
\usepackage{psfrag}
\usepackage{graphicx}
\usepackage{pinlabel}
\usepackage{graphicx}
\usepackage{amsmath}
\usepackage{amssymb}
\usepackage{amscd}
\usepackage{autobreak}
\usepackage{picinpar}
\usepackage{times}
\usepackage{pb-diagram}
\usepackage{graphicx}
\usepackage{wrapfig}
\usepackage{xspace}
\usepackage{hyperref}
\usepackage{url}
\usepackage{caption}
\usepackage{subcaption}
\usepackage{fancyhdr}
\usepackage{overpic}
\usepackage{color}
\usepackage[square,comma,sort&compress,numbers]{natbib}
\usepackage{latexsym}
\usepackage{mathrsfs}
\usepackage{amsfonts}
\usepackage{hyperref}
\usepackage{amsthm}
\usepackage{appendix}
\usepackage{pdfsync}

\theoremstyle{definition}
\newtheorem{theorem}{Theorem}[section]
\newtheorem{lemma}[theorem]{Lemma}
\newtheorem{proposition}[theorem]{Proposition}
\newtheorem{corollary}[theorem]{Corollary}
\newtheorem{definition}[theorem]{Definition}

\newtheorem{remark}[theorem]{Remark}
\newtheorem*{theorem*}{Theorem}
\def\qed{\hfill{Q.E.D.}\smallskip}

\begin{document}

\title{\bf Rigidity and deformation of generalized sphere packings on 3-dimensional manifolds with boundary}
\author{Xu Xu, Chao Zheng}

\date{\today}

\address{School of Mathematics and Statistics, Wuhan University, Wuhan, 430072, P.R. China}
 \email{xuxu2@whu.edu.cn}

\address{School of Mathematics and Statistics, Wuhan University, Wuhan 430072, P.R. China}
\email{czheng@whu.edu.cn}


\begin{abstract}
Motivated by Guo-Luo's generalized circle packings on surfaces with boundary \cite{GL2},
we introduce the generalized sphere packings on 3-dimensional manifolds with boundary.
Then we investigate the rigidity of the generalized sphere packing metrics.
We prove that the generalized sphere packing metric is determined by the combinatorial scalar curvature.
To find the hyper-ideal polyhedral metrics on 3-dimensional manifolds with prescribed combinatorial scalar curvature,
we introduce the combinatorial Ricci flow and combinatorial Calabi flow for the generalized sphere packings on 3-dimensional manifolds with boundary.
Then we study the longtime existence and convergence for the solutions of these combinatorial curvature flows.
\end{abstract}

\maketitle

\section{Introduction}

In his investigation of hyperbolic metrics on 3-dimensional manifolds, Thurston (\cite{Thurston}, Chapter 13) introduced the circle packings with non-obtuse  intersection angles on closed surfaces and proved Andreev-Thurston theorem,
which generalizes Andreev's work on the sphere \cite{Andreev 1,Andreev 2} and Koebe's work \cite{Koebe} for circle packings on the sphere.
To study the 3-dimensional analogy of the circle packings on closed surfaces, Cooper-Rivin \cite{CR} introduced the sphere packings on a 3-dimensional closed manifold, which was proved to be locally rigid \cite{CR,Glickenstein 1,Glickenstein 2,Rivin} and globally rigid \cite{Xu JDG}.
Motivated by Thurston's work \cite{Thurston},
Guo-Luo \cite{GL2} introduced the generalized circle packings on surfaces with boundary using different kinds of hyperbolic cosine laws.
This motivates us to study the generalized sphere packings on 3-dimensional manifolds with boundary.

Suppose $\Sigma$ is a compact 3-dimensional manifolds with boundary $\partial\Sigma$ consisting of $N$ connected components.
By coning off each boundary component of $\Sigma$ to be a point,
one can obtain a compact 3-dimensional space, denoted by $\widetilde{\Sigma}$.
There are exactly $N$ cone points $\{v_1,...,v_N\}$ in $\widetilde{\Sigma}$ and $\widetilde{\Sigma}-\{v_1,...,v_N\}$ is homeomorphic to $\Sigma-\partial\Sigma$.
An ideal triangulation $\mathcal{T}$ of $\Sigma$ is a triangulation $\widetilde{\mathcal{T}}$ of $\widetilde{\Sigma}$ such that the vertices of the triangulation are exactly the cone points $\{v_1,...,v_N\}$.
By Moise \cite{Moise}, every compact 3-dimensional manifold $\Sigma$ can be ideally triangulated.
Denote $st(v_1,...,v_N)$ as the open star of the vertices $\{v_1,...,v_N\}$ in the second barycentric subdivision of the triangulation $\widetilde{\mathcal{T}}$.
Then $\Sigma$ is homeomorphic to $\widetilde{\Sigma}-st(v_1,...,v_N)$.
Replacing each tetrahedron in $\mathcal{T}$ by a hyper-ideal tetrahedron and replacing the affine gluing homeomorphisms by isometries preserving the corresponding hyperbolic hexagonal faces,
we obtain a hyper-ideal polyhedral metric on $(\Sigma,\mathcal{T})$.
We call $(\Sigma,\mathcal{T})$ an ideally triangulated compact 3-dimensional manifold with boundary.

\begin{figure}[!ht]
\centering
\begin{overpic}[scale=1]{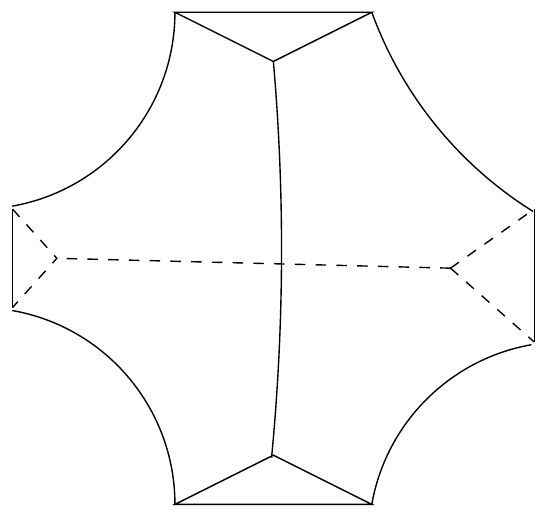}
\put (100,42) {$\triangle_h$}
\put (48,-5) {$\triangle_k$}
\put (-10,45) {$\triangle_j$}
\put (60,96) {$\triangle_i$}
\put (48,97) {$x^i_{jh}$}
\put (56,80) {$x^i_{kh}$}
\put (37,80) {$x^i_{jk}$}
\put (18,72) {$l_{ij}$}
\put (55,60) {$l_{ik}$}
\put (16,20) {$l_{jk}$}
\put (80,73) {$l_{ih}$}
\put (37,40) {$l_{jh}$}
\put (81,18) {$l_{kh}$}
\end{overpic}
\caption{Hyper-ideal tetrahedron $\sigma=\{ijkh\}$}
\label{Figure 2}
\end{figure}
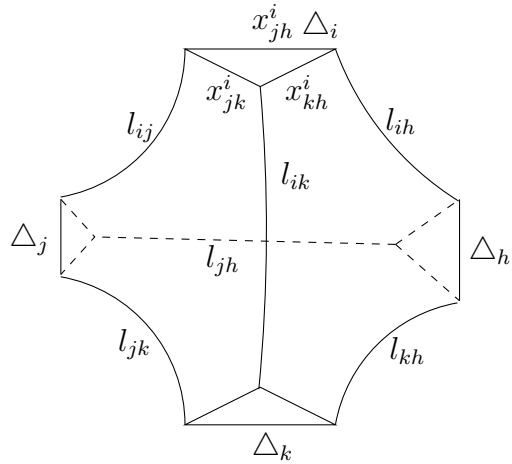

A hyper-ideal tetrahedron is a compact convex polyhedron in $\mathbb{H}^3$ that is diffeomorphic to a truncated tetrahedron in $\mathbb{E}^3$,
which has four right-angled hyperbolic hexagonal faces and four hyperbolic triangular faces.
Any triangular face is required to be orthogonal to its three adjacent hexagonal faces.
The four triangular faces isometric to hyperbolic triangles are called vertex triangles.
An edge in a hyper-ideal tetrahedron is the intersection of two hexagonal faces, and a vertex edge is the intersection of a hexagonal face and a vertex triangle.
We use $E$, $F$ and $T$ to represent the set of edges, faces and hyper-ideal tetrahedra in the triangulation $\mathcal{T}$ respectively.
A hyper-ideal tetrahedron with four vertex triangles $\triangle_\nu, \nu=i,j,k,h$ is denoted by $\sigma=\{ijkh\}$.
Using the Klein model of $\mathbb{H}^3$,
the hyper-ideal tetrahedron $\sigma=\{ijkh\}$ corresponds to a tetrahedron with four vertices $i,j,k,h$.
We call these vertices as hyper-ideal vertices of the hyper-ideal tetrahedron $\sigma=\{ijkh\}$, which corresponds to vertex triangles $\triangle_\nu, \nu=i,j,k,h$.
The set of all hyper-ideal vertices in the hyper-ideal tetrahedra is denoted by $V$.
The edge joining $\triangle_i$ to $\triangle_j$ is denoted by $\{ij\}\in E$.
The hexagonal face adjacent to $\{ij\}$, $\{jk\}$ and $\{ik\}$ is denoted by $\{ijk\}$.
The length of the edge $\{ij\}\in E$ is denoted by $l_{ij}$.
The length of the vertex edge $\triangle_i\bigcap\{ijk\}$ is denoted by $x^i_{jk}$.
Please refer to Figure \ref{Figure 2}.

\begin{definition}\label{Def: GSP metric}
Suppose $(\Sigma,\mathcal{T})$ is an ideally triangulated compact 3-dimensional manifold with boundary.
The generalized sphere packing metric is defined to be a map $r: V\rightarrow (0,+\infty)$ such that
(1) the length of the edge $\{ij\}\in E$ between two vertex triangles $\triangle_i$ and $\triangle_j$ is $l_{ij}=r_i+r_j$ and
(2) the lengths $l_{ij}, l_{ik}, l_{ih}, l_{jk}, l_{jh}, l_{kh}$ determine a non-degenerate hyper-ideal tetrahedron.
\end{definition}

The condition (2) in Definition \ref{Def: GSP metric} is called the non-degenerate condition.
Intuitively, a degenerate hyper-ideal tetrahedron is a collapsed octagon, which was called a flat hyper-ideal tetrahedron by Luo-Yang \cite{L-Y}.
Furthermore, Luo-Yang \cite{L-Y} showed that a hyper-ideal tetrahedron is non-degenerate if and only if four vertex triangles are non-degenerate, i.e., triangle inequalities hold for any vertex triangle.
Please see Section \ref{section 2} or \cite{L-Y} for more details.
Gluing these hyper-ideal tetrahedra along the faces by isometries may produce singularities at the hyper-ideal vertices on $(\Sigma,\mathcal{T})$,
which could be described by the following combinatorial scalar curvature.

\begin{definition}\label{Def: CSC}
Suppose $(\Sigma,\mathcal{T})$ is an ideally triangulated compact 3-dimensional manifold with boundary.
The combinatorial scalar curvature at a hyper-ideal vertex $i\in V$ for the generalized sphere packing metrics is defined to be
\begin{equation*}
K_i=\sum_{j;j\sim i}(2\pi-\sum_{\{ijkh\}\in T}\beta_{ij,kh})
=2\pi\chi(\Sigma_i)+\mathrm{Area}(\Sigma_i),
\end{equation*}
where $K_{ij}=2\pi-\sum_{\{ijkh\}\in T}\beta_{ij,kh}$ is called the combinatorial Ricci curvature along the edge $\{ij\}\in E$, $\beta_{ij,kh}$ is the dihedral angle at the edge $\{ij\}$ and $\Sigma_i$ is a connected component of the boundary $\partial \Sigma$.
\end{definition}

The combinatorial scalar curvature was first introduced by Cooper-Rivin \cite{CR} to study the sphere packing metrics on 3-dimensional closed manifolds. For a triangulated 3-dimensional closed manifold, the link of a vertex is a triangulated sphere,
while the link of a hyper-ideal vertex in an ideally triangulated compact 3-manifold with boundary is a triangulated surface with genus possibly bigger than 0.
This gives rise to a different form of the combinatorial scalar curvature  in Definition \ref{Def: CSC}.
Note that if the link of the  hyper-ideal vertex is a triangulated 2-sphere, then the combinatorial scalar curvature in Definition \ref{Def: CSC}
is reduced to Cooper-Rivin's combinatorial scalar curvature.

We have the following global rigidity for the generalized sphere packing metrics.

\begin{theorem}\label{Thm: result 1}
Suppose $(\Sigma,\mathcal{T})$ is an ideally triangulated compact 3-dimensional manifold with boundary.
The generalized sphere packing metric on $(\Sigma,\mathcal{T})$ is determined by its combinatorial scalar curvature.
\end{theorem}

Motivated by Chow-Luo's combinatorial Ricci flow on surface \cite{Chow-Luo}, Luo's combinatorial Ricci flow on compact 3-dimensional manifolds with boundary \cite{Luo} and Ge's combinatorial Calabi flow for Thurston's circle packing metrics on surfaces \cite{Ge1, Ge2, GX},
we introduce the following combinatorial Ricci flow and combinatorial Calabi flow for the generalized sphere packings on 3-dimensional manifolds with boundary.

\begin{definition}\label{Def: CCF}
Suppose $(\Sigma,\mathcal{T})$ is an ideally triangulated compact 3-dimensional manifold with boundary and $\overline{K}\in \mathbb{R}^N$ is a given function defined on $V$.
The combinatorial Ricci flow for the generalized sphere packing metrics is defined to be
\begin{eqnarray}\label{Eq: CRF}
\begin{cases}
\frac{dr_i}{dt}=K_i-\overline{K}_i, \\
r(0)=r_{0}.
\end{cases}
\end{eqnarray}
The combinatorial Calabi flow for the generalized sphere packing metrics is defined to be
\begin{eqnarray}\label{Eq: CCF}
\begin{cases}
\frac{dr_i}{dt}=-\Delta(K-\overline{K})_i,\\
r(0)=r_{0}.
\end{cases}
\end{eqnarray}
where $\Delta=(\frac{\partial K_i}{\partial r_j})_{N\times N}$ is the generalized discrete Laplace operator.
\end{definition}

We have the following results on the longtime existence and convergence of the combinatorial Ricci flow (\ref{Eq: CRF}) and the combinatorial Calabi flow (\ref{Eq: CCF}).

\begin{theorem}\label{Thm: result 2}
Suppose $(\Sigma,\mathcal{T})$ is an ideally triangulated compact 3-dimensional manifold with boundary and $\overline{K}\in \mathbb{R}^N$ is a given function defined on $V$.
\begin{description}
\item[(a)] If there exists a non-degenerate generalized sphere packing metric $\overline{r}$ with combinatorial scalar curvature $\overline{K}$, then the solution $r(t)$ of the combinatorial Ricci flow (\ref{Eq: CRF}) is uniformly bounded from above in $\mathbb{R}_{>0}^N$. In other words, there exists a positive constant $M$ such that $r_i(t)\leq M$ for all $i\in V$.
\item[(b)] Furthermore, for any constant $c\in (0,M)$, there exist constants $C_1=C_1(\chi(\Sigma_i),d_i,M)$ and $C_2=C_2(\chi(\Sigma_i),d_i,M,c)$ such that if $\overline{K}_i\in(C_1,C_2]$ for all $i\in V$, then the solution of the combinatorial Ricci flow (\ref{Eq: CRF}) exists for all time and converges exponentially fast to $\overline{r}$, where $\Sigma_i$ is a connected component of the boundary $\partial \Sigma$ and $d_i$ is the degree at the hyper-ideal vertex $i\in V$.
\end{description}
\end{theorem}

\begin{theorem}\label{Thm: result 3}
Suppose $(\Sigma,\mathcal{T})$ is an ideally triangulated compact 3-dimensional manifold with boundary and $\overline{K}\in \mathbb{R}^N$ is a given function defined on $V$.
If the solution of the combinatorial Calabi flow (\ref{Eq: CCF}) converges, then there exists a non-degenerate sphere packing metric $\overline{r}$ with combinatorial scalar curvature $\overline{K}$.
Furthermore, if there is a non-degenerate generalized sphere packing metric $\overline{r}$ with combinatorial scalar curvature $\overline{K}$, there exists a constant $\delta>0$ such that if $||K(r(0))-K(\overline{r})||<\delta$, then the solution of the combinatorial Calabi flow (\ref{Eq: CCF}) exists for all time and converges exponentially fast to $\overline{r}$.
\end{theorem}

The paper is organized as follows.
In Section \ref{section 2}, we prove the simply connectedness of the admissible space of the generalized sphere packing metrics for a hyper-ideal tetrahedron.
In Section \ref{section 3}, we prove the negative definiteness of the Jacobian of the combinatorial scalar curvature with respect to the generalized sphere packing metrics.
In Section \ref{section 4}, we prove the global rigidity of the generalized sphere packing metrics, i.e., Theorem \ref{Thm: result 1}.
In Section \ref{section 5}, we prove a generalization of Theorem \ref{Thm: result 3}.
By a key estimate, we also prove Theorem \ref{Thm: result 2} in this section.
\\
\\
\textbf{Acknowledgements}\\[8pt]
The author X. Xu thanks Professor Tian Yang at Texas A\&M University for helpful communications.

\section{Admissible space of generalized sphere packing metrics on a hyper-ideal tetrahedron}\label{section 2}

Suppose $\sigma=\{ijkh\}$ is a hyper-ideal tetrahedron with four vertex triangles $\triangle_\nu, \nu=i,j,k,h$ as shown in Figure \ref{Figure 2}.
The admissible space $\Omega_{ijkh}$ of the generalized sphere packing metrics on a hyper-ideal tetrahedron is the set of $(r_i, r_j, r_k, r_h)\in \mathbb{R}^4_{>0}$ such that the hyper-ideal tetrahedron with edge lengths $l_{ij}$,\ $l_{ik}$,\ $l_{ih}$,\ $l_{jk}$,\ $l_{jh}$,\ $l_{kh}$ is non-degenerate.

A basic fact in hyperbolic geometry is that given any three positive numbers, there exists a unique right-angled hyperbolic hexagon up to hyperbolic isometry with the lengths of three non-adjacent edges given by these three positive numbers \cite{Ratcliffe}.
As $l_{ij}=r_i+r_j$, then the right-angled hyperbolic hexagons on the faces always exists.
For a non-degenerate hyper-ideal tetrahedron, Luo-Yang \cite{L-Y} further gave the following proposition.
\begin{proposition}(\cite{L-Y}, Proposition 4.4)\label{Prop: L-Y}
The truncated tetrahedron $\sigma=\{ijkh\}$ with edge lengths $(l_{ij}, l_{ik}, l_{ih}, l_{jk}, l_{jh}, l_{kh})\in \mathbb{R}^6_{>0}$ determines a non-degenerate hyper-ideal tetrahedron if and only if the four vertex triangles $\triangle_\nu, \nu=i,j,k,h$ are non-degenerate, i.e., triangle inequalities hold for any vertex triangle $\triangle_\nu$.
\end{proposition}
We consider the vertex triangle $\triangle_i$ with edge lengths $x^i_{jk}, x^i_{jh}, x^i_{kh}$.
To simplify the notations, we set
\begin{equation*}
\begin{aligned}
&\sinh l_{ij}=s_{ij},\ \cosh l_{ij}=c_{ij},\ \tanh r_i=t_i.
\end{aligned}
\end{equation*}
By the cosine law for a right-angled hyperbolic hexagon, we have
\begin{equation}\label{Eq: F11}
\begin{aligned}
\cosh x^i_{jk}=& \frac{\cosh l_{ij}\cosh l_{ik}+\cosh l_{jk}}{\sinh l_{ij}\sinh l_{ik}}=\frac{c_{ij}c_{ik}+c_{jk}}{s_{ij}s_{ik}},\\
\cosh x^i_{jh}=& \frac{\cosh l_{ij}\cosh l_{ih}+\cosh l_{jh}}{\sinh l_{ij}\sinh l_{ih}}=\frac{c_{ij}c_{ih}+c_{jh}}{s_{ij}s_{ih}},\\
\cosh x^i_{kh}=& \frac{\cosh l_{ik}\cosh l_{ih}+\cosh l_{kh}}{\sinh l_{ik}\sinh l_{ih}}=\frac{c_{ik}c_{ih}+c_{kh}}{s_{ik}s_{ih}}.
\end{aligned}
\end{equation}
Note that $x^i_{jk}, x^i_{jh}, x^i_{kh}$ satisfy the triangle inequalities if and only if
\begin{equation*}
\sinh \frac{x^i_{jk}+x^i_{jh}+x^i_{kh}}{2}\sinh \frac{x^i_{jk}+x^i_{jh}-x^i_{kh}}{2}\sinh \frac{x^i_{jk}-x^i_{jh}+x^i_{kh}}{2}\sinh \frac{-x^i_{jk}+x^i_{jh}+x^i_{kh}}{2}>0.
\end{equation*}
By direct calculations, we have
\begin{alignat*}{4}
&4\sinh \frac{x^i_{jk}+x^i_{jh}+x^i_{kh}}{2}\sinh \frac{x^i_{jk}+x^i_{jh}-x^i_{kh}}{2}\sinh \frac{x^i_{jk}-x^i_{jh}+x^i_{kh}}{2}\sinh \frac{-x^i_{jk}+x^i_{jh}+x^i_{kh}}{2}\\
=&1+2\cosh x^i_{jk}\cosh x^i_{jh}\cosh x^i_{kh}-\cosh^2 x^i_{jk}-\cosh^2x^i_{jh}-\cosh^2x^i_{kh}\\
=&\frac{1}{s^2_{ij}s^2_{ik}s^2_{ih}}
[s^2_{ij}s^2_{ik}s^2_{ih}+2(c_{ij}c_{ik}+c_{jk})(c_{ij}c_{ih}+c_{jh})(c_{ik}c_{ih}+c_{kh})\\
&\ \ \ \ \ \ \ \ \ \ \ \ \ \ \ -(c_{ij}c_{ik}+c_{jk})^2s_{ih}^2-(c_{ij}c_{ih}+c_{jh})^2s_{ik}^2-(c_{ik}c_{ih}+c_{kh})^2s_{ij}^2],
\end{alignat*}
where (\ref{Eq: F11}) is used in the last line.
Furthermore, using the formula $s^2_{ij}=c_{ij}^2-1$ gives
\begin{equation}\label{Eq: Q_1}
\begin{aligned}
Q_1:=&s^2_{ij}s^2_{ik}s^2_{ih}+2(c_{ij}c_{ik}+c_{jk})(c_{ij}c_{ih}+c_{jh})(c_{ik}c_{ih}+c_{kh})\\
   &-(c_{ij}c_{ik}+c_{jk})^2s_{ih}^2-(c_{ij}c_{ih}+c_{jh})^2s_{ik}^2-(c_{ik}c_{ih}+c_{kh})^2s_{ij}^2\\
=&(c_{ij}^2-1)(c_{ik}^2-1)(c_{ih}^2-1)+2(c_{ij}c_{ik}+c_{jk})(c_{ij}c_{ih}+c_{jh})(c_{ik}c_{ih}+c_{kh})\\
 &-(c_{ij}c_{ik}+c_{jk})^2(c_{ih}^2-1)-(c_{ij}c_{ih}+c_{jh})^2(c_{ik}^2-1)-(c_{ik}c_{ih}+c_{kh})^2(c_{ij}^2-1)\\
=&c_{ij}^2+c_{ik}^2+c_{ih}^2+c_{jk}^2+c_{jh}^2+c_{kh}^2
  -c_{ij}^2c_{kh}^2-c_{ik}^2c_{jh}^2-c_{ih}^2c_{jk}^2\\
 &+2c_{ij}c_{ik}c_{jk}+2c_{ik}c_{ih}c_{kh}+2c_{jk}c_{jh}c_{kh}+2c_{ij}c_{ih}c_{jh}\\
 &+2c_{ik}c_{ih}c_{jk}c_{jh}+2c_{ij}c_{ik}c_{jh}c_{kh}+2c_{ij}c_{ih}c_{jk}c_{kh}-1.
\end{aligned}
\end{equation}
Substituting $l_{ij}=r_i+r_j$ into (\ref{Eq: Q_1}) gives
\begin{equation}\label{Eq: Q 2}
\begin{aligned}
Q_2=&\frac{1}{4\cosh^2 r_i\cosh^2 r_j\cosh^2 r_k\cosh^2 r_h}Q_1\\
=&-t^2_i-t^2_j-t^2_k-t^2_h
+2t_i2t_j+2t_i2t_k+2t_i2t_h+2t_j2t_k+2t_j2t_h+2t_k2t_h+4\\
=&(t_i+t_j+t_k+t_h)^2-2(t^2_i+t^2_j+t^2_k+t^2_h)+4.
\end{aligned}
\end{equation}

Note that $Q_2$ is symmetric in $i, j, k, h$
and hence the triangle inequalities hold for $\triangle_j, \triangle_k, \triangle_h$ if and only if $Q_2>0$.
Since $t_\nu=\tanh r_\nu\in (0,1)$ for $\nu\in \{i,j,k,h\}$, then it is direct to check that $Q_2>0$ is satisfied for any $r=(r_i,r_j,r_k,r_h)\in \mathbb{R}_{>0}^4$.
Therefore, we conclude the following result.
\begin{theorem}\label{Thm: AS connected}
The admissible space $\Omega_{ijkh}$ of the generalized sphere packing metrics on $\sigma=\{ijkh\}$ is $\mathbb{R}_{>0}^4$.
\end{theorem}

\begin{corollary}\label{Cor: AS connected}
The admissible space $\Omega=\bigcap_{\sigma\in T}\Omega_{ijkh}$ on the triangulated 3-dimensional manifold with boundary $(\Sigma,\mathcal{T})$ is $\mathbb{R}_{>0}^N$.
\end{corollary}

\section{Negative definite Jacobian for generalized sphere packing metrics}\label{section 3}
In this section, we prove the negative definiteness of the Jacobian of the combinatorial scalar curvature with respect to the generalized sphere packing metrics.
For simplicity, we set
\begin{equation*}
\begin{aligned}
\sinh r_i=s_i,\ \cosh r_i=c_i,\
\sinh(r_i+r_j+r_k)=s_{ijk},\
\cosh(r_i+r_j+r_k)=c_{ijk},\ \\
\lambda_1=t_it_j+t_it_k+t_jt_k+1,\
\lambda_2=t_it_j+t_it_h+t_jt_h+1,\
\lambda_3=t_it_k+t_it_h+t_kt_h+1.
\end{aligned}
\end{equation*}

\begin{proposition}
For the vertex triangle $\triangle_i$ with edge lengths $x^i_{jk}, x^i_{jh}, x^i_{kh}$ in a hyper-ideal tetrahedron, the dihedral angle $\beta_{ij,kh}$ at the edge $\{ij\}\in E$ satisfies
\begin{equation}\label{Eq: cos beta_ij}
\cos\beta_{ij,kh}
=\frac{c_ic_j \sqrt{c_kc_h}}{4\sqrt{c_{ijk}c_{ijh}}}
\left[Q_2-(t_i+t_j)^2+(t_k-t_h)^2\right],
\end{equation}
where $Q_2$ is defined by (\ref{Eq: Q 2}).
\end{proposition}
\proof
Using the hyperbolic cosine law for the vertex triangle $\triangle_i$ gives
\begin{equation}\label{Eq: F12}
\begin{aligned}
\cos\beta_{ij,kh}
=&\frac{-\cosh x^i_{kh}+\cosh x^i_{jk}\cosh x^i_{jh}}{\sinh x^i_{jk}\sinh x^i_{jh}}\\
=&\frac{c_{ik}c_{ih}+c_{jk}c_{jh}+c_{ij}c_{ik}c_{jh}
+c_{ij}c_{ih}c_{jk}-s^2_{ij}c_{kh}}
{\sqrt{2c_{ij}c_{ik}c_{jk}+c^2_{ij}+c^2_{ik}+c^2_{jk}-1}
\sqrt{2c_{ij}c_{ih}c_{jh}+c^2_{ij}+c^2_{ih}+c^2_{jh}-1}},
\end{aligned}
\end{equation}
where (\ref{Eq: F11}) is used in the last line.
One can refer to Lemma 4.3 in \cite{L-Y} for a detailed proof of (\ref{Eq: F12}).
Since $l_{ij}=r_i+r_j$, then
\begin{equation*}
\begin{aligned}
&2c_{ij}c_{ik}c_{jk}+c^2_{ij}+c^2_{ik}+c^2_{jk}-1\\
=&4\cosh\frac{l_{ij}+l_{ik}+l_{jk}}{2}
\cosh\frac{l_{ij}+l_{ik}-l_{jk}}{2}
\cosh\frac{l_{ij}-l_{ik}+l_{jk}}{2}
\cosh\frac{-l_{ij}+l_{ik}+l_{jk}}{2}\\
=&4c_{ijk}c_ic_jc_k.
\end{aligned}
\end{equation*}
Similarly, we have $2c_{ij}c_{ih}c_{jh}+c^2_{ij}+c^2_{ih}+c^2_{jh}-1
=4c_{ijh}c_ic_jc_h$.
The numerator of the second line of (\ref{Eq: F12}) can be calculated as follows
\begin{equation*}
\begin{aligned}
&c_{ik}c_{ih}+c_{jk}c_{jh}+c_{ij}c_{ik}c_{jh}
+c_{ij}c_{ih}c_{jk}-s^2_{ij}c_{kh}\\
=&c_{ik}(c_{ih}+c_{ij}c_{jh})+c_{jk}(c_{jh}+c_{ij}c_{ih})-s^2_{ij}c_{kh}\\
=&(c_ic_k+s_is_k)(c_ic_h+c_ic^2_jc_h+c_ic_js_js_h+c_jc_hs_is_j+c^2_js_is_h)\\
&+(c_jc_k+s_js_k)(c_jc_h+c^2_ic_jc_h+c_ic_js_is_h+c_ic_hs_is_j+c^2_is_js_h)\\
&-(c_kc_h+s_ks_h)(c_j^2s_i^2+c_i^2s_j^2+2c_ic_js_is_j)\\
=&c_i^2c_j^2c_kc_h\bigg[\frac{2(c_j^2-s_j^2)}{c^2_j}
+\frac{2(c_i^2-s_i^2)}{c^2_i}+\frac{2s_js_h}{c_jc_h}
+\frac{2s_is_k}{c_ic_k}+\frac{2s_is_h}{c_ic_h}+\frac{2s_js_k}{c_jc_k}\bigg]\\
=&c_i^2c_j^2c_kc_h(4-2t_j^2-2t_i^2+2t_jt_h+2t_it_k+2t_it_h+2t_jt_k)\\
=&c_i^2c_j^2c_kc_h[Q_2-(t_i+t_j)^2+(t_k-t_h)^2].
\end{aligned}
\end{equation*}
This completes the proof.
\qed

By (\ref{Eq: cos beta_ij}), we have
\begin{equation}\label{Eq: F13}
\sin \beta_{ij,kh}
=\sqrt{1-\cos^2 \beta_{ij,kh}}
=\frac{\sqrt{Q_3}}{4\sqrt{c_{ijk}c_{ijh}}},
\end{equation}
where
\begin{equation*}
Q_3=16c_{ijk}c_{ijh}-c^2_ic^2_jc_kc_h
[Q_2-(t_i+t_j)^2+(t_k-t_h)^2]^2.
\end{equation*}
Note that
\begin{equation*}
\begin{aligned}
&[Q_2-(t_i+t_j)^2+(t_k-t_h)^2]^2+4(t_i+t_j)^2Q_2\\
=&[Q_2+(t_i+t_j)^2+(t_k-t_h)^2]^2-4(t_i+t_j)^2(t_k-t_h)^2\\
=&(4t_it_j+2t_it_k+2t_it_h+2t_jt_k+2t_jt_h+4)^2-[2(t_i+t_j)(t_k-t_h)]^2\\
=&\frac{16}{c^2_ic^2_jc_kc_h}c_{ijk}c_{ijh},
\end{aligned}
\end{equation*}
where $c_{ijk}=\lambda_1c_ic_jc_k$ and $c_{ijh}=\lambda_2c_ic_jc_h$ are used in the last line.
This implies
\begin{equation*}
Q_3=4(t_i+t_j)^2c^2_ic^2_jc_kc_h Q_2.
\end{equation*}
Thus (\ref{Eq: F13}) can be rewritten as
\begin{equation}\label{Eq: sin beta ij}
\sin \beta_{ij,kh}
=\frac{(t_i+t_j)c_ic_j\sqrt{c_kc_h}\sqrt{Q_2}}
{2\sqrt{c_{ijk}c_{ijh}}}.
\end{equation}

\begin{proposition}
For the vertex triangle $\triangle_i$ with edge lengths $x^i_{jk}, x^i_{jh}, x^i_{kh}$ in a hyper-ideal
tetrahedron,
\begin{equation}\label{Eq: variation 1}
\frac{\partial \beta_{ij,kh}}{\partial r_k}
=-\frac{c_ic_j}{2\sqrt{Q_2}c_{ijk}c_k}
(t_i+t_j)(t_i+t_j+t_k-t_h).
\end{equation}
\end{proposition}
\proof
Differentiating (\ref{Eq: cos beta_ij}) with respect to $r_k$ gives
\begin{equation*}
\begin{aligned}
-\sin \beta_{ij,kh}\cdot \frac{\partial \beta_{ij,kh}}{\partial r_k}
=&[Q_2-(t_i+t_j)^2+(t_k-t_h)^2]\cdot\frac{\partial}{\partial r_k}\bigg(\frac{c_ic_j\sqrt{c_kc_h}}{4\sqrt{c_{ijk}c_{ijh}}}\bigg)\\
&+\frac{c_ic_j\sqrt{c_kc_h}}{4\sqrt{c_{ijk}c_{ijh}}}\cdot\left(\frac{\partial Q_2}{\partial r_k}+2(t_k-t_h)\frac{\partial t_k}{\partial r_k}\right)\\
=&\mathrm{I}+\mathrm{II}.
\end{aligned}
\end{equation*}
By direct calculations, we have
\begin{equation*}
\begin{aligned}
\mathrm{I}
=&\frac{Q_2-(t_i+t_j)^2+(t_k-t_h)^2}{16c_{ijk}c_{ijh}}
c_ic_j\sqrt{c_h}(\cdot\frac{s_k}{2\sqrt{c_k}}\cdot4\sqrt{c_{ijk}c_{ijh}}
-\sqrt{c_k}\cdot4\sqrt{c_{ijh}}\cdot\frac{s_{ijk}}{2\sqrt{c_{ijk}}})\\
=&\frac{Q_2-(t_i+t_j)^2+(t_k-t_h)^2}
{8c_{ijk}\sqrt{c_{ijk}c_{ijh}}\sqrt{c_k}}\cdot
c_ic_j\sqrt{c_h}(s_kc_{ijk}-c_ks_{ijk})\\
=&-\frac{[Q_2-(t_i+t_j)^2+(t_k-t_h)^2]c^2_ic^2_j\sqrt{c_h}
(t_i+t_j)}{8c_{ijk}\sqrt{c_{ijk}c_{ijh}}\sqrt{c_k}}
\end{aligned}
\end{equation*}
and
\begin{equation*}
\mathrm{II}
=\frac{c_ic_j\sqrt{c_kc_h}}{4\sqrt{c_{ijk}c_{ijh}}}
\left[2(t_i+t_j-t_k+t_h)\cdot\frac{\partial t_k}{\partial r_k}+2(t_k-t_h)\frac{\partial t_k}{\partial r_k}\right]
=\frac{c_ic_j\sqrt{c_h}(t_i+t_j)}
{2\sqrt{c_{ijk}c_{ijh}}c_k\sqrt{c_k}}.
\end{equation*}
Hence,
\begin{equation*}
\begin{aligned}
\frac{\partial \beta_{ij,kh}}{\partial r_k}
=&-\frac{1}{\sin \beta_{ij,kh}}(\mathrm{I}+\mathrm{II})\\
=&\frac{1}{4\sqrt{Q_2}c_{ijk}c^2_k}
\bigg(c_ic_jc_k[Q_2-(t_i+t_j)^2+(t_k-t_h)^2]-4c_{ijk}\bigg)\\
=&\frac{c_ic_j}{4\sqrt{Q_2}c_{ijk}c_k}\cdot
(-2t_i^2-2t_j^2-4t_it_j-2t_it_k+2t_it_h-2t_jt_k+2t_jt_h)\\
=&-\frac{c_ic_j}{2\sqrt{Q_2}c_{ijk}c_k}(t_i+t_j)(t_i+t_j+t_k-t_h).
\end{aligned}
\end{equation*}
\qed

\begin{proposition}
For the vertex triangle $\triangle_i$ with edge lengths $x^i_{jk}, x^i_{jh}, x^i_{kh}$ in a hyper-ideal
tetrahedron,
\begin{equation}\label{Eq: variation 2}
\frac{\partial \beta_{ij,kh}}{\partial r_i}
=\frac{c^2_jc_kc_h}{2\sqrt{Q_2}c_{ijk}c_{ijh}}\cdot H_j,
\end{equation}
where
\begin{equation}\label{Eq: H j}
\begin{aligned}
H_j
=&2t_i^2t_j^2+t_j^2t_k^2+t_j^2t_h^2+2t_i^2t_jt_k
+2t_i^2t_jt_h+2t_i^2t_kt_h+t_it^2_jt_k+t_it^2_jt_h\\
&+4t_it_jt_kt_h-2t_it_j^3-t_j^3t_k-t_j^3t_h-2t_j^2-t_k^2-t_h^2\\
&+6t_it_j+3t_it_k+3t_it_h+3t_jt_k+3t_jt_h+2t_kt_h+4.
\end{aligned}
\end{equation}
\end{proposition}
\proof
Differentiating (\ref{Eq: cos beta_ij}) with respect to $r_i$ gives
\begin{equation*}
\begin{aligned}
-\sin \beta_{ij,kh}\cdot \frac{\partial \beta_{ij,kh}}{\partial r_i}
=&[Q_2-(t_i+t_j)^2+(t_k-t_h)^2]\cdot\frac{\partial}{\partial r_i}(\frac{c_ic_j\sqrt{c_kc_h}}{4\sqrt{c_{ijk}c_{ijh}}})\\
&+\frac{c_ic_j\sqrt{c_kc_h}}{4\sqrt{c_{ijk}c_{ijh}}}\cdot\left(\frac{\partial Q_2}{\partial r_i}-2(t_i+t_j)\frac{\partial t_i}{\partial r_i}\right)\\
=&\mathrm{I_1}+\mathrm{I_2}.
\end{aligned}
\end{equation*}
By direct calculations, we have
\begin{equation*}
\begin{aligned}
\mathrm{I}_1
=&\frac{Q_2-(t_i+t_j)^2+(t_k-t_h)^2}{16c_{ijk}c_{ijh}}\cdot c_j\sqrt{c_kc_h}\\
&\times[s_i\cdot4\sqrt{c_{ijk}c_{ijh}}-4c_j
(\sqrt{c_{ijh}}\frac{s_{ijk}}{2\sqrt{c_{ijk}}}
+\sqrt{c_{ijk}}\frac{s_{ijh}}{2\sqrt{c_{ijh}}})]\\
=&\frac{c_j\sqrt{c_kc_h}[Q_2-(t_i+t_j)^2+(t_k-t_h)^2]}{8c_{ijk}c_{ijh}\sqrt{c_{ijk}c_{ijh}}}
\cdot(2s_ic_{ijk}c_{ijh}-c_is_{ijk}c_{ijh}-c_ic_{ijk}s_{ijh})\\
=&\frac{c_j\sqrt{c_kc_h}[Q_2-(t_i+t_j)^2+(t_k-t_h)^2]}{8c_{ijk}c_{ijh}\sqrt{c_{ijk}c_{ijh}}}
\cdot[-c_{ijk}\sinh(r_j+r_h)-c_{ijh}\sinh(r_j+r_k)]\\
=&-\frac{c_ic^3_jc_kc_h\sqrt{c_kc_h}}{8c_{ijk}c_{ijh}\sqrt{c_{ijk}c_{ijh}}}
\cdot[Q_2-(t_i+t_j)^2+(t_k-t_h)^2]\cdot[(t_j+t_h)\lambda_1+(t_j+t_k)\lambda_2]
\end{aligned}
\end{equation*}
and
\begin{equation*}
\begin{aligned}
\mathrm{I}_2
=&\frac{c_ic_j\sqrt{c_kc_h}}{4\sqrt{c_{ijk}c_{ijh}}}\cdot
(-4t_i+2t_k+2t_h)\cdot\frac{\partial t_i}{\partial r_i}\\
=&-\frac{c_ic^3_jc_kc_h\sqrt{c_kc_h}}{8c_{ijk}c_{ijh}\sqrt{c_{ijk}c_{ijh}}}
\cdot2(4t_i-2t_k-2t_h)\lambda_1\lambda_2.
\end{aligned}
\end{equation*}
Note that
\begin{equation*}
\begin{aligned}
&2(4t_i-2t_k-2t_h)\lambda_1\lambda_2+[Q_2-(t_i+t_j)^2+(t_k-t_h)^2][(t_j+t_h)\lambda_1+(t_j+t_k)\lambda_2]\\
=&\lambda_1[(4t_i-2t_k-2t_h)\lambda_2+(t_j+t_h)(-2t_i^2-2t_j^2+2t_it_k+2t_it_h+2t_jt_k+2t_jt_h+4)]\\
&+\lambda_2[(4t_i-2t_k-2t_h)\lambda_1+(t_j+t_k)(-2t_i^2-2t_j^2+2t_it_k+2t_it_h+2t_jt_k+2t_jt_h+4)]\\
=&\lambda_1(-2t_j^3+2t_i^2t_j+2t_i^2t_h+2t_j^2t_k+4t_it_jt_h+4t_i+4t_j-2t_k+2t_h)\\
&+\lambda_2(-2t_j^3+2t_i^2t_j+2t_i^2t_k+2t_j^2t_h+4t_it_jt_k+4t_i+4t_j+2t_k-2t_h)\\
=&2(t_i+t_j)[t_j(t_i-t_j)(\lambda_1+\lambda_2)+t_i(t_k\lambda_2+t_h\lambda_1)+2(\lambda_1+\lambda_2)]\\
&+2t_j(t_jt_k\lambda_1+t_jt_h\lambda_2+t_it_k\lambda_2+t_it_h\lambda_1)-2(t_i+t_j)(t_k-t_h)^2\\
=&2(t_i+t_j)[t_j(t_i-t_j)(\lambda_1+\lambda_2)+t_i(t_k\lambda_2+t_h\lambda_1)+2(\lambda_1+\lambda_2)-(t_k-t_h)^2]\\
&+2t_j(t_i+t_j)[t_k+t_h+t_jt_k^2+t_jt_h^2+t_it_jt_k+t_it_jt_h+2t_it_kt_h]\\
=&2(t_i+t_j) H_j.
\end{aligned}
\end{equation*}
Hence,
\begin{equation*}
\begin{aligned}
\frac{\partial \beta_{ij,kh}}{\partial r_i}
=&-\frac{1}{\sin \beta_{ij,kh}}(\mathrm{I}_1+\mathrm{I}_2)\\
=&-\frac{2\sqrt{c_{ijk}c_{ijh}}}{(t_i+t_j)c_ic_j\sqrt{c_kc_h}\sqrt{Q_2}}\cdot
\frac{-c_ic^3_jc_kc_h\sqrt{c_kc_h}}{8c_{ijk}c_{ijh}\sqrt{c_{ijk}c_{ijh}}}\cdot
2(t_i+t_j)\cdot H_j\\
=&\frac{c^2_jc_kc_h}{2c_{ijk}c_{ijh}\sqrt{Q_2}}\cdot H_j.
\end{aligned}
\end{equation*}
\qed

\begin{remark}\label{Rmk: 1}
Since $t_\nu=\tanh r_\nu\in(0,1)$ for $\nu\in\{i,j,k,h\}$, then it is easy to check $H_j>0$, which implies $\frac{\partial \beta_{ij,kh}}{\partial r_i}>0$ by (\ref{Eq: variation 2}).
Since $\beta_{ij,kh}=\beta_{ji,kh}$, then
$\frac{\partial \beta_{ij,kh}}{\partial r_i}=\frac{\partial \beta_{ji,kh}}{\partial r_i}>0$.
\end{remark}

\begin{lemma}
For the vertex triangle $\triangle_i$ with edge lengths $x^i_{jk}, x^i_{jh}, x^i_{kh}$ in a hyper-ideal
tetrahedron,
\begin{equation}\label{Eq: variation 3}
\frac{\partial \mathrm{Area}(\triangle_i)}{\partial r_j}
=-\frac{c_kc_h}{\sqrt{Q_2}c_{ijk}c_{ijh}}
[2-(t_k-t_h)^2+t_i(t_j+t_k+t_h)+t_j(t_i+t_k+t_h)].
\end{equation}
\end{lemma}
\proof
By (\ref{Eq: variation 1}), we have
\begin{equation*}
\begin{aligned}
\frac{\partial \beta_{ik,jh}}{\partial r_j}
=&-\frac{c_i^2c_kc_h}{2\sqrt{Q_2}c_{ijk}c_{ijh}}(t_i+t_k)(t_i+t_j+t_k-t_h)\lambda_2,\\
\frac{\partial \beta_{ih,jk}}{\partial r_j}
=&-\frac{c_i^2c_kc_h}{2\sqrt{Q_2}c_{ijk}c_{ijh}}(t_i+t_h)(t_i+t_j-t_k+t_h)\lambda_1.
\end{aligned}
\end{equation*}
Similar to (\ref{Eq: variation 2}), we can obtain
\begin{equation*}
\frac{\partial \beta_{ij,kh}}{\partial r_j}
=\frac{c^2_ic_kc_h}{2\sqrt{Q_2}c_{ijk}c_{ijh}}\cdot H_i,
\end{equation*}
where
\begin{equation*}
\begin{aligned}
H_i
=&2t_i^2t_j^2+t_i^2t_k^2+t_i^2t_h^2+t_i^2t_jt_k
+t_i^2t_jt_h+2t_it^2_jt_k+2t_it^2_jt_h+2t_j^2t_kt_h\\
&+4t_it_jt_kt_h-2t_i^3t_j-t_i^3t_k-t_i^3t_h-2t_i^2-t_k^2-t_h^2\\
&+6t_it_j+3t_it_k+3t_it_h+3t_jt_k+3t_jt_h+2t_kt_h+4.
\end{aligned}
\end{equation*}
Since
\begin{equation*}
\begin{aligned}
&(t_i+t_k)(t_i+t_j+t_k-t_h)\lambda_2+(t_i+t_h)(t_i+t_j-t_k+t_h)\lambda_1-H_i\\
=&4t^3_it_j+4t^2_it_kt_h+4t^2_i-4t_it_j-4t_kt_h-4\\
&+2t^3_it_h+2t^2_it_jt_h-2t_i^2t_h^2-2t_it_h-2t_jt_h+2t_h^2\\
&+2t^3_it_k+2t^2_it_jt_k-2t_i^2t_k^2-2t_it_k-2t_jt_k+2t_k^2\\
=&2(t_i^2-1)[2-(t_k-t_h)^2+t_i(t_j+t_k+t_h)+t_j(t_i+t_k+t_h)],
\end{aligned}
\end{equation*}
then
\begin{equation*}
\begin{aligned}
\frac{\partial \mathrm{Area}(\triangle_i)}{\partial r_j}
=&-(\frac{\partial \beta_{ij,kh}}{\partial r_j}
+\frac{\partial \beta_{ik,jh}}{\partial r_j}
+\frac{\partial \beta_{ih,jk}}{\partial r_j})\\
=&\frac{c^2_ic_kc_h}{2\sqrt{Q_2}c_{ijk}c_{ijh}}\cdot
2(t_i^2-1)[2-(t_k-t_h)^2+t_i(t_j+t_k+t_h)+t_j(t_i+t_k+t_h)]\\
=&-\frac{c_kc_h}{\sqrt{Q_2}c_{ijk}c_{ijh}}
[2-(t_k-t_h)^2+t_i(t_j+t_k+t_h)+t_j(t_i+t_k+t_h)].
\end{aligned}
\end{equation*}
\qed

\begin{remark}\label{Rmk: 2}
Since $t_\nu=\tanh r_\nu\in(0,1)$ for $\nu\in\{i,j,k,h\}$, then it is easy to check
\begin{equation*}
2-(t_k-t_h)^2+t_i(t_j+t_k+t_h)+t_j(t_i+t_k+t_h)>0,
\end{equation*}
which implies $\frac{\partial \mathrm{Area}(\triangle_i)}{\partial r_j}
=\frac{\partial \mathrm{Area}(\triangle_j)}{\partial r_i}<0$ by (\ref{Eq: variation 3}).
\end{remark}

\begin{theorem}\label{Thm: negative matrix}
For the vertex triangle $\triangle_i$ with edge lengths $x^i_{jk}, x^i_{jh}, x^i_{kh}$ in a hyper-ideal
tetrahedron, the Jacobian
\begin{equation*}
\Lambda_\sigma=\frac{\partial (\mathrm{Area}(\triangle_i),\mathrm{Area}(\triangle_j),\mathrm{Area}(\triangle_k),\mathrm{Area}(\triangle_h))}
{\partial(r_i,r_j,r_k,r_h)}
\end{equation*}
is symmetric and negative definite on the admissible space $\Omega_{ijkh}$.
\end{theorem}
\proof
The symmetry follows from (\ref{Eq: variation 3}).
We just need to prove the negative definiteness.
Set
\begin{equation*}
\begin{aligned}
O_j&=2-(t_k-t_h)^2+t_i(t_j+t_k+t_h)+t_j(t_i+t_k+t_h),\\
O_k&=2-(t_j-t_h)^2+t_i(t_j+t_k+t_h)+t_k(t_i+t_j+t_h),\\
O_h&=2-(t_j-t_k)^2+t_i(t_j+t_k+t_h)+t_h(t_i+t_j+t_k).
\end{aligned}
\end{equation*}
By (\ref{Eq: variation 3}), we have
\begin{align*}
\frac{\partial \mathrm{Area}(\triangle_i)}{\partial r_j}
=&-\frac{c_kc_h}{\sqrt{Q_2}c_{ijk}c_{ijh}}\cdot O_j
=-\frac{c_ic^2_jc^2_kc^2_h}{2\sqrt{Q_2}c_{ijk}c_{ijh}c_{ikh}}\cdot \frac{2O_j \lambda_3}{c^2_j},\\
\frac{\partial \mathrm{Area}(\triangle_i)}{\partial r_k}
=&-\frac{c_jc_h}{\sqrt{Q_2}c_{ijk}c_{ikh}}\cdot O_k
=-\frac{c_ic^2_jc^2_kc^2_h}{2\sqrt{Q_2}c_{ijk}c_{ijh}c_{ikh}}\cdot \frac{2O_k \lambda_2}{c^2_k},\\
\frac{\partial \mathrm{Area}(\triangle_i)}{\partial r_h}
=&-\frac{c_jc_k}{\sqrt{Q_2}c_{ijh}c_{ikh}}\cdot O_h
=-\frac{c_ic^2_jc^2_kc^2_h}{2\sqrt{Q_2}c_{ijk}c_{ijh}c_{ikh}}\cdot \frac{2O_h \lambda_1}{c^2_h}.
\end{align*}
By (\ref{Eq: variation 2}), we have
\begin{equation*}
\begin{aligned}
\frac{\partial \mathrm{Area}(\triangle_i)}{\partial r_i}
=&-(\frac{\partial \beta_{ij,kh}}{\partial r_i}
+\frac{\partial \beta_{ik,jh}}{\partial r_i}
+\frac{\partial \beta_{ih,jk}}{\partial r_i})\\
=&-\frac{c^2_jc_kc_h}{2c_{ijk}c_{ijh}\sqrt{Q_2}}\cdot H_j
-\frac{c_jc^2_kc_h}{2c_{ijk}c_{ikh}\sqrt{Q_2}}\cdot H_k
-\frac{c_jc_kc^2_h}{2c_{ikh}c_{ijh}\sqrt{Q_2}}\cdot H_h\\
=&-\frac{c_ic^2_jc^2_kc^2_h}{2\sqrt{Q_2}c_{ijk}c_{ijh}c_{ikh}}\cdot (H_j\lambda_3+H_k\lambda_2+H_h\lambda_1),
\end{aligned}
\end{equation*}
where $H_j$ is defined by (\ref{Eq: H j}),
\begin{equation*}
\begin{aligned}
H_k
=&2t_i^2t_k^2+t_j^2t_k^2+t_k^2t_h^2+2t_i^2t_jt_k
+2t_i^2t_jt_h+2t_i^2t_kt_h+t_it_jt^2_k+t_it^2_kt_h\\
&+4t_it_jt_kt_h-2t_it_k^3-t_jt^3_k-t_k^3t_h-t_j^2-2t_k^2-t_h^2\\
&+3t_it_j+6t_it_k+3t_it_h+3t_jt_k+2t_jt_h+3t_kt_h+4
\end{aligned}
\end{equation*}
and
\begin{equation*}
\begin{aligned}
H_h
=&2t_i^2t_h^2+t_j^2t_h^2+t_k^2t_h^2+2t_i^2t_jt_k
+2t_i^2t_jt_h+2t_i^2t_kt_h+t_it_jt^2_h+t_it_kt^2_h\\
&+4t_it_jt_kt_h-2t_it_h^3-t_jt^3_h-t_kt^3_h-t_j^2-t_k^2-2t_h^2\\
&+3t_it_j+3t_it_k+6t_it_h+2t_jt_k+3t_jt_h+3t_kt_h+4.
\end{aligned}
\end{equation*}
Then $\frac{\partial \mathrm{Area}(\triangle_i)}{\partial r_i}<0$ by $H_j>0, H_k>0, H_h>0$.
By Remark \ref{Rmk: 2}, we have $\frac{\partial \mathrm{Area}(\triangle_i)}{\partial r_j}<0$, $\frac{\partial \mathrm{Area}(\triangle_i)}{\partial r_k}<0$ and $\frac{\partial \mathrm{Area}(\triangle_i)}{\partial r_h}<0$.

We claim
\begin{equation}\label{Eq: F14}
\frac{\partial \mathrm{Area}(\triangle_i)}{\partial r_j}
+\frac{\partial \mathrm{Area}(\triangle_i)}{\partial r_k}
+\frac{\partial \mathrm{Area}(\triangle_i)}{\partial r_h}
-\frac{\partial \mathrm{Area}(\triangle_i)}{\partial r_i}>0.
\end{equation}
Since (\ref{Eq: F14}) is equivalent to
\begin{equation*}
-H_j\lambda_3-H_k\lambda_2-H_h\lambda_1
+\frac{2O_j \lambda_3}{c^2_j}
+\frac{2O_k \lambda_2}{c^2_k}
+\frac{2O_h \lambda_1}{c^2_h}<0,
\end{equation*}
we just need to prove
\begin{equation*}
\begin{aligned}
-H_j+\frac{2O_j}{c^2_j}&=-H_j+2O_j(1-t^2_j)<0,\\
-H_k+\frac{2O_k}{c^2_k}&=-H_k+2O_k(1-t^2_k)<0,\\
-H_h+\frac{2O_h}{c^2_h}&=-H_h+2O_h(1-t^2_h)<0.
\end{aligned}
\end{equation*}
Direct calculations give
\begin{equation*}
\begin{aligned}
-H_j+2O_j(1-t^2_j)
=&-2t_i^2t_j^2+t_j^2t_k^2+t_j^2t_h^2-2t_i^2t_jt_k
-2t_i^2t_jt_h-2t_i^2t_kt_h\\
&-3t_it^2_jt_k-3t_it^2_jt_h-4t^2_jt_kt_h
-2t_it_j^3-t_j^3t_k-t_j^3t_h-4t_it_jt_kt_h\\
&-2t_j^2-t_k^2-t_h^2-2t_it_j-t_it_k-t_it_h-t_jt_k-t_jt_h+2t_kt_h.
\end{aligned}
\end{equation*}
Since $-t_k^2-t_h^2+2t_kt_h=-(t_k+t_h)^2<0$,\
$-t_jt_k+t_j^2t_k^2=-t_jt_k(1-t_jt_k)<0$ and
$-t_jt_h+t_j^2t_h^2=-t_jt_h(1-t_jt_h)<0$,
then $-H_j+2O_j(1-t^2_j)<0$.
Similarly, we have $-H_k+2O_k(1-t^2_k)<0$ and $-H_h+2O_h(1-t^2_h)<0$.
This completes the claim.
The negative definiteness of matrix $\Lambda$ follows from the following lemma.
\begin{lemma}(\cite{Chow-Luo}, Lemma 3.10)\label{Lem: CL}
Suppose $A=[a_{ij}]_{n\times n}$ is a symmetric matrix.
If $a_{ii}>\sum_{j\neq i}|a_{ij}|$ for all indices $i$, then $A$ is positive definite.
\end{lemma}
\qed


As a consequence, we have the following result.
\begin{corollary}\label{Cor: matrix 1}
The Jacobian $\Lambda=\frac{\partial (K_i,..., K_N)}{\partial(r_i,...,r_N)}$ is symmetric and negative definite on the admissible space $\Omega$.
\end{corollary}

\section{Proof of Theorem \ref{Thm: result 1}}\label{section 4}
Corollary \ref{Cor: AS connected} and Corollary \ref{Cor: matrix 1} imply that the following energy function
\begin{equation*}
G=-\int_0^r\sum_{i=1}^N K_idr_i
\end{equation*}
is a well-defined smooth function on the admissible space $\Omega$ with $\nabla_r G=-K$.
Furthermore,
\begin{equation*}
\mathrm{Hess}_r G=-\Lambda>0,
\end{equation*}
which implies the function $G$ on $(\Sigma,\mathcal{T})$ is a strictly convex function by Corollary \ref{Cor: matrix 1}. Then the global rigidity follows from the following lemma.
\begin{lemma}\label{Lem: analysis}
If $f:\Omega \rightarrow \mathbb{R}$ is a $C^1$-smooth  strictly convex function on an open convex set $\Omega \subset \mathbb{R}^n$,
then its gradient $\nabla f:\Omega \rightarrow \mathbb{R}^n$ is injective.
Furthermore, $\nabla f$ is a smooth embedding.
\end{lemma}
\qed

\section{Combinatorial curvature flows}\label{section 5}
%

The combinatorial Ricci flow (\ref{Eq: CRF}) and the combinatorial Calabi flow (\ref{Eq: CCF}) are ODE systems with smooth coefficients.
Therefore, the solutions always exist locally around the initial time $t=0$.
We further have the following result on the longtime existence and convergence for the solutions of these combinatorial curvature flows, which includes Theorem \ref{Thm: result 3} as special cases.

\begin{theorem}\label{Thm: local convergent}
Suppose $(\Sigma,\mathcal{T})$ is an ideally triangulated compact 3-dimensional manifold with boundary and $\overline{K}\in \mathbb{R}^N$ is a function defined on $V$.
If the solution of the combinatorial Ricci flow (\ref{Eq: CRF}) or the combinatorial Calabi flow (\ref{Eq: CCF}) converges,
there exists a non-degenerate generalized sphere packing metric $\overline{r}$ with combinatorial scalar curvature $\overline{K}$.
Furthermore, if there is a non-degenerate sphere packing metric $\overline{r}$ with combinatorial scalar curvature $\overline{K}$, there exists a constant $\delta>0$ such that if $||K(r(0))-K(\overline{r})||<\delta$, then the solutions of the combinatorial Ricci flow (\ref{Eq: CRF}) and the combinatorial Calabi flow (\ref{Eq: CCF}) exist for all time and converge exponentially fast to $\overline{r}$ respectively.
\end{theorem}
\proof
Suppose $r(t)$ is a solution of the combinatorial Ricci flow (\ref{Eq: CRF}).
If $\overline{r}:=r(+\infty)=\lim_{t\rightarrow +\infty}r(t)$ exists in $\Omega$,
then $K(\overline{r})=\lim_{t\rightarrow +\infty}K(r(t))$ exists by the $C^1$-smoothness of $K$.
Furthermore, there exists a sequence $\xi_n\in(n,n+1)$ such that
\begin{equation*}
r_i(n+1)-r_i(n)=r'_i(\xi_n)=K_i(r(\xi_n))-\overline{K}_i\rightarrow 0,\ \text{as}\  n\rightarrow +\infty,
\end{equation*}
which implies $K_i(\overline{r})=\lim_{n\rightarrow +\infty}K_i(r(\xi_n))=\overline{K}_i$ for all $i\in V$.
Similarly, if the solution $r(t)$ of the combinatorial Calabi flow (\ref{Eq: CCF}) converges to $\overline{r}$ with combinatorial scalar curvature $\overline{K}$ as $t\rightarrow +\infty$, then
$K(\overline{r})=\lim_{t\rightarrow +\infty}K(r(t))$ exists by the $C^1$-smoothness of $K$.
Furthermore, there exists a sequence $\xi_n\in(n,n+1)$ such that as $n\rightarrow +\infty$,
\begin{equation*}
r_i(n+1)-r_i(n)=r'_i(\xi_n)
=\Delta(K(r(\xi_n))-\overline{K})_i\rightarrow 0,
\end{equation*}
which implies $K(\overline{r})=\overline{K}$ by the negative definiteness of the generalized discrete Laplace operator $\Delta$ in Corollary \ref{Cor: matrix 1}.

For the combinatorial Ricci flow (\ref{Eq: CRF}), set $\Gamma(r)=K-\overline{K}$.
Then $D\Gamma|_{r=\overline{r}}=\Lambda$ is negative definite by Corollary \ref{Cor: matrix 1}.
Therefore, $\overline{r}$ is a local attractor of the combinatorial Ricci flow (\ref{Eq: CRF}). Then the conclusion follows from Lyapunov Stability Theorem (\cite{Pontryagin}, Chapter 5).
Similarly, for the combinatorial Calabi flow (\ref{Eq: CCF}), set $\Gamma(r)=-\Delta(K-\overline{K})$.
Then $D\Gamma|_{r=\overline{r}}=-\Lambda^2$ is negative definite by Corollary \ref{Cor: matrix 1}.
Hence, $\overline{r}$ is a local attractor of the combinatorial Calabi flow (\ref{Eq: CCF}).
The conclusion follows from Lyapunov Stability Theorem (\cite{Pontryagin}, Chapter 5).
\qed

Theorem \ref{Thm: local convergent} gives the longtime existence and convergence of the solutions of the combinatorial Ricci flow (\ref{Eq: CRF}) and the combinatorial Calabi flow (\ref{Eq: CCF}) for initial value with small energy.
For general initial value, we need the following lemma.

\begin{lemma}\label{Lem: key estimate}
Let $\triangle_i$ be the vertex triangle with edge lengths $x^i_{jk}, x^i_{jh}, x^i_{kh}$ in a hyper-ideal tetrahedron.
Suppose the generalized sphere packing metric $r$ is uniformly bounded from above in $\mathbb{R}_{>0}^4$, i.e., $r_\nu\leq M$ for any $\nu\in\{i,j,k,h\}$.
Then there exists a positive number $\mathcal{C}_1(M)$ such that $\mathrm{Area}(\triangle_i)> \mathcal{C}_1(M)$.
Furthermore, if $r_i<c$ for any $c\in(0,M)$, then there exists a positive number $\mathcal{C}_2(M,c)$ such that  $\mathrm{Area}(\triangle_i)> \mathcal{C}_2(M,c)>\mathcal{C}_1(M)$.
\end{lemma}
\proof
By the proof of Theorem \ref{Thm: negative matrix},
we have $\mathrm{Area}(\triangle_i)$ is decreasing in $r_i,r_j,r_k,r_h$.
Thus for $r_\nu\in(0,M],\ \forall\nu\in\{i,j,k,h\}$, we have
\begin{equation*}
\min\mathrm{Area}(\triangle_i)
=\mathrm{Area}(\triangle_i)|_{(r_i,r_j,r_k,r_h)=(M,M,M,M)}.
\end{equation*}
The formula (\ref{Eq: cos beta_ij}) can be rewritten as
\begin{equation}\label{Eq: F17}
\cos\beta_{ij,kh}(r_i,r_j,r_k,r_h)
=\frac{2-t_j^2-t_i^2+t_it_k+t_it_h+t_jt_k+t_jt_h}
{2\sqrt{(t_it_j+t_it_k+t_jt_k+1)(t_it_j+t_it_h+t_jt_h+1)}},
\end{equation}
which implies
\begin{equation*}
\cos\beta_{ij,kh}|_{(r_i,r_j,r_k,r_h)=(M,M,M,M)}=\frac{1+\tanh^2 M}{1+3\tanh^2 M}:=\widetilde{C}_1(M).
\end{equation*}
Similarly,
$\cos\beta_{ik,jh}|_{(r_i,r_j,r_k,r_h)=(M,M,M,M)}
=\widetilde{C}_1(M)$ and
$\cos\beta_{ih,jk}|_{(r_i,r_j,r_k,r_h)=(M,M,M,M)}
=\widetilde{C}_1(M)$.
Set $\mathcal{C}_1(M)=\pi-3\mathrm{arccos}(\widetilde{C}_1(M))$.
Then
\begin{align*}
\mathrm{Area}(\triangle_i)
\geq&\mathrm{Area}(\triangle_i)|_{(r_i,r_j,r_k,r_h)=(M,M,M,M)}\\
=&\pi-(\beta_{ij,kh}+\beta_{ik,jh}+\beta_{ih,jk})|_{(r_i,r_j,r_k,r_h)=(M,M,M,M)}\\
=&\mathcal{C}_1(M).
\end{align*}
Furthermore, if $r_i<c$ for any $c\in(0,M)$, then
$\mathrm{Area}(\triangle_i)
>\mathrm{Area}(\triangle_i)|_{(r_i,r_j,r_k,r_h)=(M,M,M,M)}$.
By (\ref{Eq: F17}), we have
$\cos\beta_{ij,kh}|_{(r_i,r_j,r_k,r_h)=(c,M,M,M)}
=\frac{1}{2}+\frac{1-\tanh^2 c}{2(1+\tanh^2 M+2\tanh c\tanh M)}:=\widetilde{C}_2(M,c)$.
Similarly,
$\cos\beta_{ik,jh}|_{(r_i,r_j,r_k,r_h)=(c,M,M,M)}
=\widetilde{C}_2(M,c)$
and
$\cos\beta_{ih,jk}|_{(r_i,r_j,r_k,r_h)=(c,M,M,M)}
=\widetilde{C}_2(M,c)$.
Set $\mathcal{C}_2(M,c)=\pi-3\mathrm{arccos}(\widetilde{C}_2(M,c))$.
Then
\begin{align*}
\mathrm{Area}(\triangle_i)
>&\mathrm{Area}(\triangle_i)|_{(r_i,r_j,r_k,r_h)=(c,M,M,M)}\\
=&\pi-(\beta_{ij,kh}+\beta_{ik,jh}+\beta_{ih,jk})|_{(r_i,r_j,r_k,r_h)=(c,M,M,M)}\\
=&\mathcal{C}_2(M,c).
\end{align*}
\qed

\noindent\textbf{Proof\ of\ Theorem\ \ref{Thm: result 2}}
Suppose there exists a non-degenerate generalized sphere packing metric $\overline{r}$ with combinatorial scalar curvature $\overline{K}$.
Corollary \ref{Cor: AS connected} and Corollary \ref{Cor: matrix 1} imply that the following energy function
\begin{equation*}
\widetilde{G}(r)=-\int_{\overline{r}}^r\sum_{i=1}^N (K_i-\overline{K}_i)dr_i.
\end{equation*}
is a well-defined smooth function on the admissible space $\Omega$.
By direct calculations, we have
\begin{equation*}
\nabla_{r_i}\widetilde{G}(r)=-(K_i-\overline{K}_i),
\end{equation*}
which implies the combinatorial Ricci flow (\ref{Eq: CRF}) is a negative gradient flow of the function $\widetilde{G}(r)$.
Then $\widetilde{G}(r)$ is a strictly convex function defined on $\mathbb{R}_{>0}^N$ because of $\mathrm{Hess}_r \widetilde{G}=-\Delta>0$.
Furthermore, $\widetilde{G}(\overline{r})=0,\ \nabla \widetilde{G}(\overline{r})=0$,
which implies $\widetilde{G}(r)\geq \widetilde{G}(\overline{r})=0$ and $\lim_{r\rightarrow \infty}\widetilde{G}(r)=+\infty$ by the convexity of $\widetilde{G}(r)$.
On can also refer to Lemma 4.6 in \cite{Ge-Xu} for a proof for this fact.
Moreover,
\begin{equation*}
\frac{d\widetilde{G}(r(t))}{dt}
=\sum_{i=1}^{N}\frac{\partial \widetilde{G}}{\partial r_i}\frac{d r_i}{dt}
=-\sum_{i=1}^{N}(K_i-\overline{K}_i)^2\leq0,
\end{equation*}
which implies that the function $\widetilde{G}(r(t))$ is decreasing along the combinatorial Ricci flow (\ref{Eq: CRF}) and $0\leq \widetilde{G}(r(t))\leq \widetilde{G}(r(0))$.
Hence, the solution $r(t)$ of the combinatorial Ricci flow (\ref{Eq: CRF}) is bounded in $\mathbb{R}_{>0}^N$.
We claim that $r_i(t)$ is uniformly bounded from below in $\mathbb{R}_{>0}$ along the combinatorial Ricci flow (\ref{Eq: CRF}).
We shall prove the theorem assuming the claim and then prove the claim.

By the claim, we have the solution of the combinatorial Ricci flow (\ref{Eq: CRF}) exists for all time and $\widetilde{G}(r(t))$ converges.
Moreover, there exists a sequence $\xi_n\in(n,n+1)$ such that as $n\rightarrow +\infty$,
\begin{equation*}
\widetilde{G}(r(n+1))-\widetilde{G}(r(n))
=(\widetilde{G}(r(t))'|_{\xi_n}=\nabla \widetilde{G}\cdot\frac{dr_i}{dt}|_{\xi_n}
=-\sum^N_{i=1}(K_i-\overline{K}_i)^2|_{\xi_n} \rightarrow 0.
\end{equation*}
Then $\lim_{n\rightarrow +\infty}K_i(r(\xi_n))=K_i
=K_i(\overline{r})$ for all $i\in V$.
By $\{r(t)\}\subset\subset \mathbb{R}_{>0}^N$, there exists $r^*\in \mathbb{R}_{>0}^N$ and
a subsequence of $\{r(\xi_n)\}$, still denoted as $\{r(\xi_n)\}$ for simplicity,
such that $\lim_{n\rightarrow \infty}r(\xi_n)=r^*$, which implies
$K_i(\overline{r})=\lim_{n\rightarrow +\infty}K_i(r(\xi_n))=K_i(r^*)$.
This further implies $\overline{r}=r^*$ by Theorem \ref{Thm: result 1}.
Therefore, $\lim_{n\rightarrow \infty}r(\xi_n)=\overline{r}$.
Set $\Gamma(r)=K-\overline{K}$, then $D\Gamma|_{r=\overline{r}}$ has $N$ negative eigenvalues, and hence $\overline{r}$ is a local attractor of the combinatorial Ricci flow (\ref{Eq: CRF}).
Then the conclusion follows from Lyapunov Stability Theorem (\cite{Pontryagin}, Chapter 5).

We use Ge-Xu's trick in \cite{Ge-Xu 17} to prove the claim.
Suppose that there exists at least one hyper-ideal vertex $i\in V$ such that $\lim_{t\rightarrow T}r_i(t)=0^+$ for $T\in(0,+\infty)$.
Since the solution $r(t)$ of the combinatorial Ricci flow (\ref{Eq: CRF}) is bounded in $\mathbb{R}_{>0}^N$,
then $r(t)$ is uniformly bounded from above in $\mathbb{R}_{>0}^N$, denoted by $M$.
By Lemma \ref{Lem: key estimate}, there exists a positive
number $\mathcal{C}_1(M)$ such that $\mathrm{Area}(\triangle_i)\geq\mathcal{C}_1(M)$,
which implies
\begin{equation*}
K_i=2\pi\chi(\Sigma_i)+\mathrm{Area}(\Sigma_i)
=2\pi\chi(\Sigma_i)+d_i\mathrm{Area}(\triangle_i)
\geq2\pi\chi(\Sigma_i)+d_i\mathcal{C}_1(M),
\end{equation*}
where $\Sigma_i$ is a connected component of the boundary $\partial \Sigma$ and $d_i$ is the degree at the hyper-ideal vertex $i\in V$.
Furthermore, if $r_i<c$ for any $c\in(0,M)$, then there exists a positive constant $\mathcal{C}_2(M,c)$ such that $\mathrm{Area}(\triangle_i)> \mathcal{C}_2(M,c)>\mathcal{C}_1(M)$,
which implies
\begin{equation*}
K_i>2\pi\chi(\Sigma_i)+d_i\mathcal{C}_2(M,c)
>2\pi\chi(\Sigma_i)+d_i\mathcal{C}_1(M).
\end{equation*}
Set $C_1(\chi(\Sigma_i),d_i,M)=2\pi\chi(\Sigma_i)+d_i\mathcal{C}_1(M)$
and
$C_2(\chi(\Sigma_i),d_i,M,c)=2\pi\chi(\Sigma_i)+d_i\mathcal{C}_2(M,c)$.
Note that $\overline{K}_i\in (C_1(\chi(\Sigma_i),d_i,M),C_2(\chi(\Sigma_i),d_i,M,c)]$.
Choose a time $t_0\in (0,T)$ such that $r_i(t_0)<c$, this can be done because $\lim_{t\rightarrow T}r_i(t)=0^+$. Set $a=\inf\{t<t_0|r_i(s)<c, \forall s\in (t,t_0]\}$, then $r_i(a)=c$.
Note that $\frac{dr_i}{dt}=K_i-\overline{K}_i>0$ on $(a,t_0]$, we have $r_i(t_0)>r_i(a)=c$, which contradicts $r_i(t_0)<c$.
Therefore, $r_i(t)$ is uniformly bounded from below in $\mathbb{R}_{>0}$.
\qed

\noindent\textbf{Data\ availability\ statements}
Data sharing not applicable to this article as no datasets were generated or analysed during the current study.


\begin{thebibliography}{50}
\setlength{\itemsep}{2pt} \small

\bibitem{Andreev 1} E. M. Andreev, \emph{On convex polyhedra of finite volume in Loba\u{c}evski\u{i} spaces}, {Math. USSR-Sb.} 12 (1970), 255-259.
\bibitem{Andreev 2} E. M. Andreev, \emph{On convex polyhedra in Loba\u{c}evski\u{i} spaces}, {Math. USSR-Sb.} 10 (1970), 412-440.


\bibitem{Chow-Luo} B. Chow, F. Luo, \emph{Combinatorial Ricci flows on surfaces}, J. Differential Geom. 63 (2003), no. 1, 97-129.

\bibitem{CR} D. Cooper, I. Rivin, \emph{Combinatorial scalar curvature and rigidity of ball packings}. Math. Res. Lett. 3 (1996), 51-60.

\bibitem{Ge1} H. Ge, \emph{Combinatorial methods and geometric equations}, Thesis (Ph.D.)-Peking University, Beijing. 2012. (In Chinese).

\bibitem{Ge2} H. Ge, \emph{Combinatorial Calabi flows on surfaces}, Trans. Amer. Math. Soc. 370 (2018), no. 2, 1377-1391.

\bibitem{GX}  H. Ge, X. Xu, \emph{$2$-dimensional combinatorial Calabi flow in hyperbolic background geometry}. Differential Geom. Appl. 47 (2016), 86-98.

\bibitem{Ge-Xu 17} H. Ge, X. Xu, \emph{A discrete Ricci flow on surfaces with hyperbolic background geometry}, Int. Math. Res. Not. IMRN 2017, no. 11, 3510-3527.

\bibitem{Ge-Xu} H. Ge, X. Xu, \emph{On a combinatorial curvature for surfaces with inversive distance circle packing metrics}. J. Funct. Anal. 275 (2018), no. 3, 523-558.

\bibitem{Glickenstein 1} D. Glickenstein, \emph{A combinatorial Yamabe flow in three dimensions}, Topology 44 (2005), No. 4, 791-808.


\bibitem{Glickenstein 2} D. Glickenstein, \emph{A maximum principle for combinatorial Yamabe flow}, Topology 44 (2005), No. 4, 809-825.

\bibitem{GL2} R. Guo, F. Luo, \emph{Rigidity of polyhedral surfaces, II}, Geom. Topol. 13 (2009), no. 3, 1265-1312.

\bibitem{Koebe} P. Koebe, \emph{Kontaktprobleme der konformen Abbildung}. Ber. S\"{a}chs. Akad. Wiss. Leipzig, Math. Phys. Kl. 88 (1936), 141-164.

\bibitem{Luo} F. Luo, \emph{A combinatorial curvature flow for compact 3-manifolds with boundary}, Electron. Res. Announc. Am. Math. Soc. 11 (2005) 12-20.

\bibitem{L-Y} F. Luo, T. Yang, \emph{Volume and rigidity of hyperbolic polyhedral 3-manifolds}, J. Topol. 11 (1) (2018) 1-29.

\bibitem{Moise} E. Moise, \emph{Affine structures in 3-manifolds V}. Ann. Math. (2) 56 (1952), 96-114.

\bibitem{Pontryagin} L.S. Pontryagin, \emph{Ordinary differential equations}, Addison-Wesley Publishing Company Inc., Reading, 1962.


\bibitem{Ratcliffe} J.G. Ratcliffe, \emph{Foundations of hyperbolic manifolds}. Second edition. Graduate Texts in Mathematics, 149, xii+779 pp. Springer, New York (2006). ISBN: 978-0387-33197-3; 0-387-33197-2.

\bibitem{Rivin} I. Rivin, \emph{An extended correction to ``Combinatorial Scalar Curvature and Rigidity of Ball Packings,''} (by D. Cooper and I. Rivin), \href{https://arxiv.org/abs/0302069v2}
    {arXiv:0302069v2[math.GT]}.

\bibitem{Thurston} W. Thurston, \emph{Geometry and topology of $3$-manifolds}, Princeton lecture notes 1976,\href{http://www.msri.org/publications/books/gt3m}
    {http://www.msri.org/publications/books/gt3m}.

\bibitem{Xu JDG} X. Xu, \emph{On the global rigidity of sphere packings on 3-dimensional manifolds}, J.Differential Geom. 115(2020), no. 1, 175-193.

\end{thebibliography}
\end{document}